\magnification\magstep 1
\documentstyle{amsppt}
\def\a{\alpha} \def\b{\beta} \def\d{\delta}
\def\g{\gamma}  \def\k{\kappa}
 \def\s{\sigma}
 \def\sm{\setminus} \def\o{\omega}
\def\to{\rightarrow}  \def\imp{\Rightarrow}
\def\r{\upharpoonright}  \def\cl{\overline}
  \def\0{\emptyset}
 
\def\<{\langle} \def\>{\rangle}
\topmatter
\nologo
\title  Spaces having a small diagonal 
\endtitle

\author Gary Gruenhage \endauthor \address Department of
Mathematics, Auburn University, AL 36849\endaddress

\email garyg\@mail.auburn.edu \endemail

\abstract We obtain several results and examples concerning the general question ``When must a
space with a small diagonal have a $G_\d$-diagonal?".  In particular, we show 
(1) every compact metrizably fibered space with 
a small diagonal is metrizable;  
(2) there are consistent  examples 
of regular Lindel\"of (even hereditarily Lindel\"of) spaces with 
a small diagonal but no $G_\delta$-diagonal; (3) every first-countable
hereditarily Lindel\"of space with a small diagonal has a $G_\d$-diagonal;
(4) assuming $CH$, every Lindel\"of $\Sigma$-space with a small diagonal has 
a countable network;
(5) whether  countably compact spaces 
with a small diagonal are metrizable depends on your set theory;
(6) there is a locally compact space with a small diagonal but no 
$G_\delta$ diagonal. 
\endabstract

\thanks Research partially supported by NSF DMS-9704849 \endthanks

\endtopmatter

\document

\head 0.  Introduction \endhead

According to M. Hu\v sek [H], a space $X$ has a {\it small diagonal}
\footnote"$^1$"{Hu\v sek actually used the more descriptive ``$\o_1$-inaccessible diagonal", 
but the term ``small diagonal", which was suggested by E. van Douwen, seems to have
become more popular.} if 
for every uncountable subset $Y$ of $X^2\setminus \Delta$, there is
an open  set $U\supset \Delta$ such that $Y\setminus U$ is 
uncountable.  Clearly a space with a $G_\delta$-diagonal has a 
small diagonal; the question is, for what classes of spaces  does small diagonal imply 
$G_\delta$-diagonal?  

This question for compact $T_2$-spaces is a 
well-known and still not completely solved problem of Hu\v sek, who proved that, 
assuming the Continuum Hypothesis ($CH$), compact spaces of countable tightness
having a small diagonal are metrizable (recall metrizability is equivalent to 
$G_\d$-diagonal for compact or even countably compact spaces--see, e.g., [Gr]).  
Dow [D] showed that this result holds in any model obtained by adding
Cohen reals over a model of $CH$, and Zhou[Z] proved, assuming $CH$
plus what he called ``Fleissner's Axiom", that compact spaces having a
small diagonal are metrizable.
Later, Juhasz and Szentmikl\"ossy [JS] showed in $ZFC$ that 
non-countably-tight compact spaces
cannot have a small diagonal (because they have convergent
$\o_1$-sequences).  This implies that the countable tightness assumption 
in the above results of Hu\v sek and Dow can be omitted, i.e., 
compact spaces having a small
diagonal are metrizable in models of  $CH$ or in their extensions by 
Cohen forcing.

Hu\v sek also asked the corresponding question for $\o_1$-compact 
and other spaces.  Zhou [Z] studied the
question for Lindel\"of, countably compact, and locally compact spaces, and obtained, under
$MA+\neg CH$, a locally compact example, and a hereditarily Lindel\"of, but non-regular,
example of a space with a small diagonal and no $G_\d$-diagonal.
Arhangel'ski and Bella [AB] generalized the afore-mentioned $CH$ result
for compact spaces 
to Lindel\"of spaces which are perfect pre-images of metrizable spaces.   

Here we show: 
(1) every compact metrizably fibered space with 
a small diagonal is metrizable;  
(2) there are consistent  examples 
of regular Lindel\"of (even hereditarily Lindel\"of) spaces with 
a small diagonal but no $G_\delta$-diagonal; (3) every first-countable
hereditarily Lindel\"of space with a small diagonal has a $G_\d$-diagonal;
(4) assuming $CH$, every Lindel\"of $\Sigma$-space with a small diagonal has 
a countable network;
(5) whether  countably compact spaces 
with a small diagonal are metrizable depends on your set theory;
(6) there is a locally compact space with a small diagonal but no 
$G_\delta$ diagonal. 
 The Lindel\"of $\Sigma$-space result (4) answers a 
question of Arhangel'skii, statement (5) answers questions of Zhou and 
Shakhmatov, and (2) and (6) answer questions left open by Zhou.  

In the sequel we mention several questions which remain open, including: 
\roster
\item"(1)" Is it true in $ZFC$ that every compact space (or Lindel\"of $\Sigma$-space)
with a small diagonal is metrizable?
\item"(2)" Is there in $ZFC$ a Lindel\"of space with a small diagonal but no $G_\d$-diagonal?
\item"(3)" Can there be a first-countable countably compact, or first-countable Lindel\"of, 
space with a small diagonal but no $G_\d$-diagonal\footnote"$^2$"{\ O.
Pavlov[P] has recently obtained 
a positive solution to the countably compact question.}?
\endroster

Unless stated otherwise, all spaces are assumed to be regular and $T_1$. 
 
\head 1.  Preliminaries\endhead

It will be helpful in the sequel to make some observations about
small diagonals that are probably well known to anyone who has
considered this property.  First, note that any open superset of
the diagonal in $X^2$ contains an open set of the form
$\bigcup_{U\in \Cal U}U^2$ for some open cover $\Cal U$ of $X$. 
The following is then easy to prove:

\proclaim{Lemma 1.1} The following are equivalent for a space $X$:
\roster
\item"(a)" $X$ has a small diagonal;
\item"(b)" Whenever $\Cal D$ is a uncountable collection of
doubletons in $X$, there is an open cover $\Cal U$ such that, for
uncountably many $d=\{d_1,d_2\}\in \Cal D$, $$\forall U\in\Cal
U(d\not\subset U).$$
\endroster
\endproclaim

We recall here that it is well-known and easy to see that a space 
with a small diagonal cannot contain a convergent $\o_1$-sequence, i.e.,
a sequence $\<x_\a\>_{\a<\o_1}$ such that every neighborhood of some
point $x$ contains all but countably many $x_\a$'s.  For then the set of
points $(x_\a,x)$, $\a<\o_1$, in $X^2\sm \Delta$ is readily seen to
witness the failure of the small diagonal property (from the definition
or from Lemma 1.1(b)).  It is also clear that the small diagonal property is hereditary.   

\proclaim{Proposition 1.2} Of the statements below, 
(c)$\imp$(b)$\imp$(a).  
If $X$ is Lindel\"of, all are equivalent.
\roster
\item"(a)" $X$ has a small diagonal;
\item"(b)"  Whenever $\{\{d_\a^0,d_\a^1\}:\a<\o_1\}$ is a
collection of  doubletons of $X$, there are disjoint closed
sets $H_0$ and $H_1$ with $d_\a^0\in H_0$ and $d_\a^1\in H_1$ for
uncountably many $\a<\o_1$;
\item"(c)"  Whenever $\{\{d_\a^0,d_\a^1\}:\a<\o_1\}$ is a
collection of  doubletons of $X$, there is a co-zero (if
$X$ is 0-dimensional we can say clopen) set $U$ such that
$d_\a^0\in U$ and $d_\a^1\not\in U$ for uncountably many $\a\in
\o_1$.  
\endroster
\endproclaim
\demo{Proof}  That (c) implies (b) is clear from the fact that co-
zero sets are $F_\sigma$.  For (b) implies (a), note that the
complement of $H_0\times H_1$ is an open superset of the diagonal.

Now assume $X$ is Lindel\"of.   We show (a) implies (c). 
Let $\Cal D=\{\{d_\a^0,d_\a^1\}:\a<\o_1\}$ be a collection of
 doubletons of $X$.  There is an open cover $\Cal U$ of $X$
satisfying the conditions of Lemma 1.1(b) with respect to $\Cal D$. 
Since $X$ is regular and Lindel\"of, hence completely regular, 
we may assume $\Cal U$ 
is countable and consists of co-zero sets (clopen sets if $X$ is
0-dimensional).  For each $\a<\o_1$, there is some $U_\a\in \Cal U$
with $d_\a^0\in U_\a$.  Note that $d_\a^1\not\in U_\a$.  Now the
result follows since $U_\a $ must be the same member of $\Cal U$
for uncountably many $\a$.  \qed
\enddemo

\head 2.  Lindel\"of $\Sigma$-spaces and compact metrizably fibered spaces \endhead

A space $X$ is a Lindel\"of $\Sigma$-space [N] if it is a continous image of a perfect 
pre-image of a separable metric space; equivalently, there is a countable collection
$\Cal F$ of closed sets and a cover $\Cal K$ by compact sets such that, whenever $U$ is 
an open superset of some $K\in \Cal K$, then $K\subset F\subset U$ for some $F\in \Cal F$.
The class of Lindel\"of $\Sigma$-spaces 
can be viewed as a common generalization of the class of compact spaces and separable
metric spaces.  Every $K$-analytic space (see, e.g., [RJ] for the definition) is 
in this class.  A Lindel\"of $\Sigma$-space has a $G_\delta$-diagonal iff
it has a countable network (and hence iff it is a continuous image of a 
separable metric space). 

As a generalization of the result for compact spaces, Arhangel'skii and Bella [AB] proved, 
assuming $CH$, that if $X$ is a perfect pre-image of a 
separable metric space,
and has a small diagonal, then $X$ is metrizable.    
(Bennett and Lutzer[BL] showed, however, that there are paracompact 
--but necessarily non-Lindel\"of -- perfect pre-images of metric spaces
 having a small diagonal but no
$G_\delta$-diagonal. ) 

Thus the following question, due to Arhangel'skii [A, Problem 70], 
is natural:  Is it true, or at
least consistent, that a Lindel\"of $\Sigma$-space $X$ with a small diagonal must have a 
countable network (equivalently, 
must be a continuous image of a separable metric space)?  The question
was repeated by Tkachuk [T$_1$], who answered it affirmatively in case
$X$ is a space of the form $C_p(Y)$, i.e., all continuous real-valued
valued functions on $Y$ with the topology of pointwise convergence.  

In this section, we solve part of Arhangel'skii's question by showing 
that the answer is positive under $CH$.  
The main result here is in fact the following 
theorem of $ZFC$, which has the $CH$ result as a corollary.

\proclaim{Theorem 2.1}  Suppose $X$ is a regular 
Lindel\"of $\Sigma$-space, witnessed by the countable collection
$\Cal F$ of closed sets and cover $\Cal K$ by compact sets.  If every member
of $\Cal K$ is metrizable,  and $X$ has a small diagonal, then $X$ has a
countable network.  
\endproclaim

Before embarking on the proof, we first establish the following useful fact.

\proclaim{Lemma 2.2} Suppose $X$ is a regular space, and 
that $\Cal F$ and $\Cal K$ satisfy the conditions
for $X$ to be a Lindel\"of $\Sigma$-space, where $\Cal F$ is closed under 
finite intersections.  Let $\Cal K^*$ be the collection 
of all non-empty intersections from the collection $\{K\cap F: K\in \Cal K,
F\in \Cal F\}$.  Then $\Cal F$ and $\Cal K^*$ also satisfy these
conditions.  
\endproclaim

\demo{Proof}  Let $\Cal H$ be the collection of all closed sets $H$ such
that for any open superset $U$ of $H$, there is some $F\in \Cal F$ with
$H\subset F\subset U$.  Note that both $\Cal K$ and $\Cal F$ are contained 
in $\Cal H$.  We need to show $\Cal K^*\subset \Cal H$.  

We first show that if $H_0,H_1\in \Cal H$, then $H_0\cap H_1\in \Cal H$. 
Suppose $H_0\cap H_1\subset U$, where $U$ is open.  Since $X$ is regular
Lindel\"of, hence normal, there are disjoint open sets $V_0,V_1$
containing $H_0\setminus U$ and $H_1\setminus U$, respectively.  
Now there are $F_0,F_1\in \Cal F$
containing $H_0,H_1$, respectively, and contained in $U\cup V_0,U\cup
V_1$, respectively.  Then $F_0\cap F_1$ contains $H_0\cap H_1$, and it
is easy to check that $F_0\cap F_1\subset U$.  

It follows that we may assume $\Cal K$ is closed under intersections
with members of $\Cal F$.  Thus it remains to check that every non-empty
intersection of members of $\Cal K$ is in $\Cal H$.  By the above
paragraph, this is true for finite intersections.  Suppose $\kappa$ is
an infinite
cardinal and every non-empty intersection of fewer than $\k$ members 
from $\Cal K$ is in $\Cal H$.  The lemma follows if we can show that 
whenever $\{K_\a:\a<\k\}\subset \Cal K$ and $\0\neq\bigcap_{\a<\k}K_\a$,  
then $\bigcap_{\a<\k}K_\a\in \Cal H$.  But this holds because any open
superset of $\bigcap_{\a<\k}K_\a$ contains $\bigcap_{\a<\b}K_\a$ for
some $\b<\a$, and $\bigcap_{\a<\b}K_\a\in \Cal H$ by the inductive 
assumption. \qed
\enddemo

\demo{Proof of Theorem 2.1} Assume $X$, $\Cal F$, and $\Cal K$ satisfy the 
hypotheses of the theorem. 

For each $K\in \Cal K$, let $\Cal F_K=\{F\in
\Cal F:K\subset F\}$.  If $p\in X$ and $\Cal N$ is a collection of sets, 
let us say that {\it $\Cal N$ generates a network at $p$} if the collection
of all finite intersections of members of $\{N\in \Cal N: p\in N\}$ is a 
network at $p$.  First we show:

{\it Claim. For each $K\in \Cal K$, there is a countable collection
$\Cal U_K$ of cozero sets such that $\Cal U_K\cup\Cal F_K$ generates a 
network at each point of $K$}.  

{\it Proof of Claim.}  Take $K\in \Cal K$.  Since $K$ is separable
metric, we may let $\Cal B_K$ be a countable (relative) 
base for the subspace $K$.  For each pair $B_0,B_1\in \Cal B$ having 
disjoint closures, since $X$ is normal we can choose a disjoint cozero 
sets in $X$ containing them.  Let $\Cal U_K$ be the collection of    
these chosen cozero sets.  Suppose $p\in K$, and 
  consider an open neighborhood $U$ of $p$. 
 Let $N_0,N_1,\dots$ list all members of 
$\Cal U_K\cup \Cal F_K$ which contain $p$.  If no finite intersection of
$N_i$'s is contained in $U$, then choose 
$x_n\in \bigcap_{i\leq n}N_i\setminus U$.  Since the $x_n$'s diagonalize 
through members of $\Cal F_K$, and $\Cal F_K$ is an outer network for 
the compact set $K$, the $x_n$'s  must have some limit point, say $q$,  
in $K$. Certainly $q\not\in U$.  Thus there are $B_0,B_1\in \Cal B_K$ 
having disjoint closures and containing $p$ and $q$, respectively.
Then there is a cozero set $V$ in $\Cal U_K$ containing $p$ whose
closure misses $q$.  But all but finitely many $x_n$'s are in $V$,
contradiction.  Thus $\Cal U_K\cup \Cal F_K$ generates a network at $p$,
which proves the claim.  

 We now  observe that to complete the proof of the theorem, it
 suffices to show that there is a 
countable collection $\Cal U$ of cozero sets such that 
$\Cal U\cup \Cal F$ separates points in the 
$T_1$ sense.  For, if such $\Cal U$ exists, then for each $U\in \Cal U$ we 
can add to $\Cal F$ a countable collection of closed sets whose union is
$U$, and close up under finite intersections.  Then  
every singleton is in the collection $\Cal K^*$ defined in Lemma 2.2, 
whence $\Cal F$ is a countable network
for $X$.

Suppose then that no such collection $\Cal U$ exists.  Pick $K_0\in \Cal
K$, and let $\Cal U_{K_0}$ be as in the Claim.  By assumption, $\Cal
U_{K_0}\cup \Cal F$ is not $T_1$-separating, so there are distinct
points $x_1,y_1$ such that every member of $\Cal
U_{K_0}\cup \Cal F$ which contains $x_1$ also contains $y_1$.  It
follows that every member of $\Cal K$ which contains $x_1$ also contains
$y_1$; in particular, there is some $K_1\in \Cal K$ with $x_1,y_1\in K_1$.  
Then let $\Cal U_{K_1}$ be as in the Lemma 2.2.  

Suppose $\a<\omega_1$ and we have defined $K_{\b}$ for all $\b<\a$, and points
$x_\b\neq y_\b\in K_\b$ for $0<\b<\a$, such that every member of 
$\Cal F\cup(\bigcup_{\g<\b}\Cal U_{K_\g})$ 
which contains $x_\b$ also contains $y_\b$.     
Since the collection 
$\Cal F\cup(\bigcup_{\g<\a}\Cal U_{K_\g})$ is not $T_1$-separating, we can
find $x_\a\neq y_\a$ such that every member of the collection which
contains $x_\a$ also contains $y_\a$.  Then choose some $K_\a\in \Cal K$
containing $x_\a$, and note that $K_\a$ must contain $y_\a$ too.  

Thus we can define $x_\a,y_\a$, and $K_\a$ as above for all $\a<\o_1$.  
Since $X$ has a small diagonal and is Lindel\"of, there are disjoint
closed sets $H_0,H_1$ and an uncountable subset $W$ of $\o_1$ such that
$x_\a\in H_0$ and $y_\a\in H_1$ for all $\a\in W$.  

Now consider $\a\in
W$.  For each $p\in K_\a$, assuming as we may that the collections 
$\Cal F$ and $\Cal U_{K_\a}$ are closed under finite intersections,
by the Claim there is some $U_p^\a\in \Cal U_{K_\a}$ and some $F_p^\a\in 
\Cal F_{K_\a}$ such that $p\in U_p^\a$ and $U_p^\a\cap F_p^\a$ misses
either $H_0$ or $H_1$.  By compactness, there are finite subcollections 
$U_0^\a,U_1^\a,\dots, U_{n_\a}^\a$ of $\Cal U_{K_\a}$ and $F_0^\a, F_1^\a,  
\dots,F_{n_\a}^\a$ of $\Cal F_{K_\a}$ such that the $U_i^\a$'s cover
$K_\a$ and each $U_i^\a\cap F_i^\a$ misses either $H_0$ or $H_1$.  
Choose $F^\a\in \Cal F_{K_\a}$ such that $K_\a\subset
F^\a\subset \bigcup_{i\leq n_\a}U_i^\a$ and $F^\a\subset \bigcap_{i\leq
n_\a}F_i^\a$.  

There are $\a_1<\a_2$ with $F^{\a_1}=F^{\a_2}=F$.  Then $K_{\a_1}\cup 
K_{\a_2}\subset F\subset \bigcup_{i\leq n_{\a_1}}U_i^{\a_1}$.  Choose 
$i\leq n_{\a_1}$ with $x_{\a_2}\in U_i^{\a_1}$.  Since $U_i^{\a_1}\cap 
F_i^{\a_1}$ misses either $H_0$ or $H_1$, and $x_{\a_2}\in H_0$, it must
miss $H_1$.  Thus $x_{\a_2}$ is in $U_i^{\a_1}\cap 
F_i^{\a_1}$ but $y_{\a_2}$ is not, contradicting the way $x_{\a_2}$ and
$y_{\a_2}$ were chosen.  This contradiction completes the proof of the
theorem.  \qed
\enddemo

Since compact spaces with a small diagonal are metrizable under $CH$, 
the following corollary is immediate.

\proclaim{ Corollary 2.3 (CH)} Every regular Lindel\"of $\Sigma$-space with a small 
diagonal has a countable network.
\endproclaim

Fremlin [F$_1$] showed that, assuming $MA+\neg CH$, if every compact
subset of a $K$-analytic space $X$ is metrizable, then $X$ is analytic.
 Fremlin's result fails under $CH$, but we have the
following corollary to our theorem.

\proclaim{ Corollary 2.4(CH)} Every regular 
$K$-analytic space with a small diagonal 
is analytic.  
\endproclaim

\demo{Proof} Let $X$ be regular $K$-analytic space with a small
diagonal.  Then $X$ is a Lindel\"of $\Sigma$-space, hence has a countable
network.  A $K$-analytic space with a countable network is analytic (see, e.g.,
[RJ; Theorem 5.5.1]).  \qed
\enddemo
 
Our theorem can also be applied to the class of compact metrizably fibered
spaces.  According to Tkachuk [T$_2$], $X$ is {\it metrizably fibered} if
there is a continous mapping $f:X\to M$ for some metrizable space $M$,
such that each point-inverse is metrizable.  The class of metrizably
fibered compacta contains the Alexandroff duplicate of the interval, the
Alexandroff double arrow space, and many variations of these spaces.
This class has been a rich source of examples in topology (see, e.g.,
[W] or [GN]).  
The following corollary, this time a $ZFC$ result, 
shows that no member of this class can provide an
answer to Hu\v sek's question about compact spaces with a small diagonal.

\proclaim{Corollary 2.5} A metrizably fibered compact space with a
small diagonal must be metrizable. 
\endproclaim

\demo{Proof} Let $X$ be compact, and let $f:X\to M$ be a continuous map
from $X$ onto the metrizable space $M$, with $f^{-1}(m)$ metrizable for 
each $m\in M$.  Then $M$ has a countable base $\Cal B$.  Let 
$\Cal F=\{f^{-1}(B):B\in \Cal B\}$ and $\Cal K=\{f^{-1}(m):m\in M\}$.  
Then $X$, $\Cal F$, and $\Cal K$ satisfy all the hypotheses of  
Theorem 2.1.  So $X$ has a countable network, hence is metrizable. \qed
\enddemo

{\bf Remark 2.6.} Some well-known members of the class of compact metrizably fibered spaces
are perfectly normal, equivalently, hereditarily Lindel\"of 
(e.g., the double arrow space).  We note
here that it is a corollary to Theorem 3.6 that in fact any 
perfectly normal compact space with a
small diagonal is metrizable (in $ZFC$).

\head 3. General Lindel\"of spaces \endhead

Zhou \cite{Z} gave an example, under $MA+\neg CH$, of a Hausdorff, non-regular,
(hereditarily) Lindel\"of space 
which has a small diagonal but no $G_\delta$-diagonal.  It has
remained unsolved 
whether or not there could be a regular example.  In this section,
we give two different 
constructions of consistent regular examples.  One exists in a model of 
$\neg CH$, and is hereditarily Lindel\"of, the other exists in some
models of $CH$ (including 
$V=L$).  The latter example 
shows a contrast with the compact case, where with $CH$ small diagonal
does imply $G_\d$-diagonal.

Our example consistent with $CH$ is obtained by modifying a
construction due to Shelah for 
building an example of a Lindel\"of space of cardinality $\o_2$
($=c^+$ since $CH$ holds) 
in which each point is a $G_\d$.  The space cannot have a
$G_\d$-diagonal, for if it did, it
would have a weaker separable metrizable topology (see, e.g., [Gr;
Corollary 2.9]) and 
hence could not have cardinality greater than $\frak c$.  We don't know if
the space as defined by 
Shelah always has a small diagonal, but we will show that it is
possible to modify the forcing 
to make sure the diagonal will be small.   Let us also remark that
an easier construction 
of a large size Lindel\"of points $G_\d$ space due to Gorelic
\cite{Go}
 does not seem to lend itself
to a similar modification.  

We will closely follow the presentation due to Juh\'asz \cite{J$_1$} of Shelah's
example.  We recall the 
following definition:

\definition{Definition} A map $f:X^2\to 2$ is call {\it flexible} 
if for any distinct $x,y\in X$ 
and $i,j\in 2$ there is $z\in X$ such that $f(z,x)=i$ and
$f(z,y)=j$.  
\enddefinition

 For any $x\in X$
and $i\in 2$, put 
$$A_x^i=\{y\in X: y\neq x\text{ and }f(x,y)=i\}.$$ Let $\tau_f^i$
be the topology on $X$ having as 
subbase sets of the form $A_x^i\cup \{x\}$ and their complements.

Let us call a map $f:\o_2^2\to 2$ {\it very flexible} if it is
flexible, and 
\medskip

\noindent ($*$)\ \ whenever $\a<\b\leq \g<\o_2$, 
there exists $\d\in (\g,\g+\o)$ with $f(\d,\a)\neq f(\d,\b)$.  
\medskip

Shelah shows that there is a countably closed $\omega_2$-c.c. poset
forcing
a flexible function $f:\o_2^2\to 2$ such that the topologies
$\tau_f^i$ are Lindel\"of
with points $G_\d$.  We will see that if $f$ is very flexible, then
the resulting space 
has a small diagonal.  We then show that it is possible to modify
the forcing so that 
$f$ is very flexible.  

\proclaim{Lemma 3.1}  If $f:\o_2^2\to 2$ is very flexible, then $\o_2$
with either of  
the topologies $\tau_f^i$  has a small diagonal.  
\endproclaim
\demo{Proof}  Suppose $\{\a_\mu, \b_\mu\}_{\b<\o_1}\subset
[\o_2]^2$.  Choose $\g<\o_2$ 
with $\g>\a_\mu+\b_\mu$ for every $\mu<\o_1$.  Since $f$ is very
flexible, for each 
$\mu<\o_1$ there is $\d_\mu\in (\g,\g+\o)$ with
$f(\d_\mu,\a_\mu)\neq f(\d_\mu,\b_\mu)$.
Thus for some $\d$,  $f(\d,\a_\mu)\neq f(\d,\b_\mu)$ for
uncountably
many $\mu$.  For these $\mu$, $A_{\d}^i\cup \{\d\}$ 
contains exactly one of $\a_\mu,\b_\mu$ and is clopen in
$\tau_f^i$.   Thus $(\o_2, \tau_f^i)$ has a small diagonal by 
Proposition 1.2.
\qed
\enddemo

We continue to follow \cite{J$_1$}.  Let $\operatorname{Fn}(I<J)$ denote
the set of finite partial functions from $I$ into $J$. 
For $s\in \operatorname{Fn}(\o_2,2)$, put $$U_s=\bigcap \{A_\a^{s(\a)}\cup
\{\a\}:\a\in
\operatorname{dom}(s)\}$$
and let $\Cal U_f =\{U_s:\s\in \operatorname{Fn}(\o_2,2)\}$.  $\Cal U_f$ 
is said to be {\it Lindel\"of} if 
every cover of $\o_2$ by members of $\Cal U_f$ has a countable
subcover.    
As is shown in \cite{J$_1$}, if $\Cal U_f$ is Lindel\"of, then the
topologies
$\tau_f^i$ are Lindel\"of with points $G_\d$.  So it remains to
prove
the following Theorem 3.2, which is precisely Theorem 1.6 of \cite{J$_1$} with
"very
flexible" in place of "flexible".   The proof is also the similar
to that given in 
\cite{J$_1$}, with one extra condition on members of the poset 
 so that $F$ will be sure to be very flexible.  However, it's not
completely 
obvious that the same proof works with this extra condition, so it
will be 
necessary to define the poset and give several key parts of the 
argument. But we will not repeat here the parts that are clearly
not affected by the extra condition.  

\proclaim{Theorem 3.2}  Con($ZF$)$\imp$ Con($ZFC+CH+\exists$ a very
flexible 
$F:\o_2^2\to 2$ for which $\Cal U_F$ is Lindel\"of). 
\endproclaim
\demo{Proof}  Assume the ground model $V$ satisfies $ZFC+CH$.   A condition
$p$ will determine a countable subset $A^p$ of $\o_2$ and a
countable fragment $f^p:(A^p)^2\to 2$ of $F$.  Now for $s\in
\operatorname{Fn}(\o_2,2)$, put $$U_s^p=\{x\in A:\forall z\in \operatorname{dom}(s)(z\neq x\imp 
f^p(z,x)=s(z))\}.$$ 

A condition $p\in P$ is a triple $p=\<A,f,T\>$ satisfying:
\roster
\item"(i)"$A\in  [\o_2]^\o$;
\item"(ii)" $f:A^2\to 2$;
\item"(iii)"$|T|\leq\o$ and $\forall B\in T(B\subset \operatorname{Fn}(A,2)\ \&\
\cup 
\{U^p_s:s\in B\}=A)$;
\item"(iv)" $\forall B\in T\forall \d\in A\forall y\in (A\sm \d)
\forall h\in \operatorname{Fn}(A\sm\d,2) \exists s\in B (y\in U_{s\r\d}^p\ \&\
s\cup h\in \operatorname{Fn}(A,2))$;
\item"(v)" whenever $\a<\b\leq \g\in A$, 
there exists $\d\in A\cap (\g,\g+\o)$ with $f(\d,\a)\neq f(\d,\b)$.

\endroster

Define a 3-place relation $E^p$ on $A$ by 
$$E^p(\d,y,z)\iff[\d\leq y,z\ \&\ \forall x\in A\cap
\d(f(x,y)=f(x,z))].$$  Then if $p'=\<A',f',T'\>\in P$, we say 
$p'\leq p$ iff $A\subset A'$, $f\subset f'$, $T\subset T'$, and
$E^p\subset E^{p'}$.   

This poset is the same as the poset $P$ given in \cite{J$_1$} except
for the additional condition (v) (and the quite trivial but
technically 
useful change in (i) disallowing $|A|<\o$).
Clearly condition (v) will in
the end give us that $F$ is very flexible, once we have shown that
everything else goes through as before.

Note that condition (v) does not affect whether or not $E^p\subset
E^{p'}$, since $E^p$ is determined by values $f(x,y)$ for $x<y$,
i.e., 
values above the diagonal, 
while condition (v) is determined by values below the diagonal.
Let us note also that condition (iv) is determined by values of $f$
above the diagonal, since $\d\leq y$ there and 
the truth of $y\in U_{s\r\d}^p$ depends 
on what $f(z,y)$ is for $z\in \s\r\d$.  Keeping this in mind will
greatly simplify our task ahead.

The proof that $P$ is $\o_1$-complete is easy and the same as in
\cite{J$_1$}.  What will require some work is showing that Lemmas 1.9
and 1.10 of \cite{J$_1$} still hold.  Juh\'asz's Lemma 1.9 is:

\proclaim{Lemma 3.3} $P$ satisfies the $\o_2$-c.c..
\endproclaim
\demo{Proof}  As in \cite{J$_1$}, by a $\Delta$-system argument, 
the proof boils down to showing that two conditions $p=\<A,f,T\>$
and $p'=\<A',f',T'\>$ are compatible whenever:
\roster
\item"(a)"$\Delta=A\cap A'<A\sm\Delta,A'\sm\Delta$;
\item"(b)"$f\r\Delta^2=f'\r\Delta^2$;
\item"(c)"$\operatorname{type}(A\sm\Delta)=\operatorname{type}(A'\sm\Delta)$;
\item"(d)"if $z\in A\sm\Delta$ and $z'\in A'\sm\Delta$ are such
that $\operatorname{type}(A\cap z)=\operatorname{type}(A'\cap z')$, then $f(x,z)=f(x,z')$ for all
$x\in \Delta$;
\item"(e)" The natural bijection $\theta:A\to A'$ induces a
bijection from 
$T$ to $T'$.
\endroster

To this end, a function $g:(A\cup A')^2\to 2$ is constructed which
extends both $f$ and $f'$, and so that $q=\<A\cup A',g,T\cup T'\>$
is in $P$ and extends both $p$ and $p'$. The function $g$ needs to
be defined on
$$[(A\sm\Delta)\times(A'\sm\Delta)\cup(A'\sm\Delta)\times(A\sm
\Delta)].$$
We can define $g$ on $(A\sm\Delta)\times(A'\sm\Delta)$ as in
\cite{J$_1$}, 
since as we noted above the extra condition (v) only depends on 
values of $g$ below the diagonal. 
Note that this will also get condition (iii) holding 
for $B\in T$.

{\it Definition of $g$ on $(A'\sm\Delta)\times(A\sm\Delta)$.}   We
are 
below the diagonal here, so we don't need to worry about condition 
(iv) or $E^{q}$.  Satisfying condition (iii) for $B'\in T'$, and 
condition (v), are our only concern.  

Enumerate in type $\o$  all 4-tuples $\<B', a, a', a''\>$ where 
$B'\in T'$, $a\in A\sm\Delta$, $a'\in A'\sm\Delta$, and 
$a''\in (A\cup A')\cap (a'+1)$.   The function $g$ is defined 
by induction, finitely many values at a time.  
Suppose at stage $n$, we are given the $n^{th}$ 4-tuple $\<B', a,
a',
a''\>$.  Let $k'(x,a)=g(x,a)$ for the finitely many $x\in
A'\sm\Delta$ 
for which $g(x,a)$ has been defined.  Let $B=\theta^{-1}(B')$, and
let
$k=k'\circ \theta \in \operatorname{Fn}(A\sm\Delta, 2)$.  Apply condition (iv) for
$p$ 
with this $B$, and with $\d=min(A\sm\Delta)$, $y=a$, and $h=k$ to
get 
$s\in B$ satisfying $$a\in U_{s\r\d}^p\text{ and }
s\cup k\in \operatorname{Fn}(A,2).$$  
Let $s'=s\circ \theta^{-1}\in T'$ and note that $s'\in B'$.
  Since $s$ and $k$ are compatible, so
are $s'$ and $k'$.  Also, $s\r \delta= s\r\Delta=s'\r\Delta$.     
Thus $a\in U^p_{s\r\d}$ implies $g(z,a)=f(z,a)=s(z)=s'(z)$ for all
$z\in
\Delta \cap \operatorname{dom}(s)$.  We can extend $g$ so that now $g(x,a)=(s'\cup
k')(x)$
for every $x\in (A'\sm\Delta)\cap \operatorname{dom}(s'\cup k')$.  This is
consistent 
(by the use of $k'$) with how $g$ was defined at previous steps of
the
induction, and it puts $a\in U_{s'}^q$.  So, when we're done, 
the conclusion of (iii) will hold for $B'$.  

To make sure (v) will hold, at this same stage look at $(a,a',a")$
and
choose some $\d\in (a',a'+\o)\cap A'$ such that $g(\d,a)$ has not
yet 
been defined.  (Since (v) holds for $p'$, $(a',a'+\o)\cap A'$ must
be 
infinite.)  Then simply make $g(\d,a)$ different from $g(\d,a")$. 

Finally, let $g(a',a)=0$ if it is not yet defined.
\medskip

This completes the definition of $g$, and the verification that $q$
satisfies conditions (iii) and (v).  Verification of the other
conditions and that $q$ extends $p$ and $p'$ is the same as in 
\cite{J$_1$}, and for the most part is 
easily observed from the fact
that above the diagonal this $g$ is the same as the $g$ there. \qed
 \enddemo

Essentially what remains to show now is the folowing analogue of 
Lemma 1.10 of \cite{J$_1$}.  The difference is that here we use
$A\cup [z,z+\o)$ instead of $A\cup \{z\}$; we can't put the latter
because condition (v) implies that $(z,z+\o)\cap A$ is infinite 
whenever $z\in A$.  As in \cite{J$_1$}, it's part (b) which implies
that the resulting generic $F$ is flexible.

\proclaim{Lemma 3.4}  Suppose $p=\<A,f,T\>\in P$.  Then
\roster
\item"(a)" for every $z\in \omega_2\sm A$ there is an extension 
$g=(A\cup [z,z+\o))^2\to 2$ of $f$ such that $q=\<A\cup
[z,z+\o),g,T\>\in P$ and $q\leq p$;
\item"(b)" If $z\in \o_2\sm (\cup A+1)$, $\d,y\in A$ with $\d\leq
y$ and
$h\in \operatorname{Fn}(A\sm \d,2)$, there is an extension $g=(A\cup
[z,z+\o))^2\to 2$ 
of $f$ such that $q=\<A\cup
[z,z+\o),g,T\>\in P$, $q\leq p$, $g(x,z)=h(x)$ for every $x\in
\operatorname{dom}(h)$
and moreover $E^q(\d,y,z)$ holds.  
\endroster
\endproclaim

\demo{Proof}  (a){\it Case 1.  $|(z,z+\o)\cap A|=\o$. }  In this
case,
follow the proof in \cite{J$_1$} to first extend to $A\cup \{z\}$,
but to
obtain (v) also index the triples $(\a,\b,\g)\in (A\cup\{z\})^3$
with
$z\in \{\a,\b\}$.  (Note that for $z=\g$ condition (v) will hold
simply 
because it holds for $p$.)  Then at the $n^{th}$ step, after
extending 
$g$ to finitely many more points as in \cite{J$_1$}, extend $g$ to 
some point $(z,y)$ where $g$ is not yet defined so that (v) will be
witnessed for the $n^{th}$ triple.  Then repeat the above for
$z_1,z_2,...$, where the $z_i$'s enumerate $(z,z+\o)\sm A$.

{\it Case 2. $|(z,z+\o)\cap A|<\o$. }  Note that, by (v),
$(z,z+\o)\cap A=\0$.
List $[z,z+\o)$ as $z_0,z_1,...$ and add them in one at a time as
in
Case 1.  Triples $(\a,\b,\g)$ with $\g\in A$ can be taken care of
as
before.  Triples with $\g\in [z,z+\o)$ can easily 
be taken care of using the
fact that, when we consider $z_n$,  $g(z_n, y)$ for $y<z_n$ can be
defined arbitrarily.  

(b) Add in $z+n$, $n=0,1,...$ one at a time, defining 
$g(x,z+n)$ as in \cite{J$_1$} (with $z=z+n$).  Define $g(z+n,y)$ for
$y<z+n$
so that in the end (v) will hold for any  $\g\in [z,z+\o)$.  
 \qed
\enddemo

The remainder of the proof of Theorem 3.2 follows as in \cite{J$_1$}, so that
completes our
argument. \qed

\enddemo

{\bf Remark.} As with Shelah's example, the above construction can be done 
using an $(\o_1,1)$  morass with built-in $\diamondsuit$ sequence, 
which exists in Godel's constructible universe $L$ (see [V], Theorem 5.3.2);
essentially what goes on is that under these assumptions there is a 
filter $G$ on the partial order $P$ meeting enough dense sets 
to obtain the desired function $F:\o_2^2\to 2$.

We now turn to the hereditarily Lindel\"of example, which is built
from 
an $HFC^2$ in $2^{\o_2}$.  Recall (see \cite{J$_2$}or \cite{J$_3$}) that 
an uncountable subset $F$ of $ 2^{\o_2}$ is $HFC^k$ ($k\in \o$) if
for every
uncountable subset $W$ of $F^k$, and for every  
collection $\{\<\s^i_n\>_{i<k}\}_{n\in\o}\subset (\operatorname{Fn}(\o_2,2))^k$
with
$\operatorname{dom}(\s^i_n)=H_n$ for all $i<k$ and $n\in\o$, where the $H_n$'s are
disjoint
and have the same cardinality,  there is some $n\in \o$ and some
$k$-tuple
$\<g^i\>_{i<k}\in W$ with $\s^i_n\subset g^i$ for all $i<k$. 
Recall also that any $HFC (=HFC^1)$ is hereditarily Lindel\"of,  
and there is an $HFC^2$ (in fact a strong $HFC$) of
cardinality
$\o_2$ in $2^{\o_2}$ in any model obtained by adding $\o_2$-many
Cohen reals.

\proclaim{Theorem 3.5}  Suppose there is an $HFC^2$ $F$ in $2^{\o_2}$
of
cardinality $\o_2$.   Then there is a hereditarily Lindel\"of
space with a small diagonal but no $G_\d$-diagonal.  
\endproclaim
\demo{Proof}  For convenience, we may index $F$ as
$\{f_\a^e:\a<\o_2,e<2\}$.
Now define $g_\a:\o_2\to\o_2$ by 
$g_\a(\b)=f_\a^0(\b)$ if $\b<\a$ and $g_\a(\b)=f_\a^1(\b)$ if
$\b\geq\a$.

Our example is the subspace
$X=\{g_\a:\a<\o_2\}\cup\{f_\a^0:\a<\o_2\}$
 of $2^{\o_2}$.  We'll
prove $X$ is $HFC$ and hence hereditarily Lindel\"of.
Supose $\{\s_n\}_{n\in\o}\subset  \operatorname{Fn}(\o_2,2)$ is such that the
$\s_n$'s have
disjoint domains of the same cardinality, and 
$W\in [\o_2]^{\o_1}$.  Then, by $F$ being $HFC^2$,  there are
$\a\in W$ and
$n\in\o$
such that $f_\a^e\supset\s_n$ for each $e<2$.  Note that this
implies that
$g_\a\supset\s_n$.  It follows that both $\{f_\a^0:\a<\o_2\}$
and $\{g_\a:\a<\o_2\}$ are $HFC$-spaces, thus $X$ is too and is
hereditarily Lindel\"of.

We show $X$ has no $G_\d$-diagonal.  Suppose on the contrary that
$(\Cal
G_n)_{n\in\o}$ is a $G_\d$-diagonal sequence for $X$.  Since $X$ is
Lindel\"of, we may assume $\Cal G_n=\{[\s]:\s\in \Sigma_n\}$, where
$\Sigma_n\in [\operatorname{Fn}(\o_2,2)]^\o$.  Choose $\a<\o_2$ with
$$\a>\bigcup\{\operatorname{dom}(\s):\s\in \cup_{n\in\o}\Sigma_n\}.$$  
Then since $g_\a\r\a=f_\a^0\r\a$, we have $g_\a\in [\s]$ whenever
$\s\in
\Sigma_n$ and $f_\a^0\in [\s]$.  Thus $g_\a\in st(f_\a^0, \Cal
G_n)$ for
all $n$, contradiction.  

It remains to show that $X$ has a small diagonal.  To this end, let
$\{h_\g^0,h_\g^1\}_{\g<\o_1}\subset [X]^2$.  There is $\mu<\o_2$
sufficiently large so that if $h_\g^e=g_\a$, then $\a<\mu$.  It
follows
that $$\{h_\g^0\r(\o_2\sm\mu), h_\g^1\r(\o_2\sm\mu)\}=
\{f_\a^i\r(\o_2\sm\mu), f_\b^{j}\r(\o_2\sm\mu)\}$$ where either 
$\a\neq\b$ or $i\neq j$.  Let  $\o_1=\bigcup_{\b<\o_1}W_\b$, where
the
$W_\b$'s are disjoint.  Applying $HFC^2$ to each $W_\b$ and
$\{\<\mu+n,i\>_{i<2}\}_{n\in\o}$, we see that there are $\a_\b\in
W_\b$ and
$n_\b\in \o$ such that  $h_{\a_\b}^0(\mu+n)=0$ and $
h_{\a_\b}^1(\mu+n)=1$. 
It follows that  for some 
$\d\in (\mu,\mu+\o)$, $h_\g^0(\d)\neq h_\g^1(\d)$ for uncountably
many $\d$.  
So $X$ has a small diagonal by Proposition 1.2. \qed
 
\enddemo

{\bf Remark.} 
Any hereditarily Lindel\"of regular space has cardinality not
greater than
$2^\o$,  so any model which contains an $HFC$ of cardinality $\o_2$ 
cannot
satisfy $CH$.  We don't know if  there is a space having the
properties of the
above example which is consistent with $CH$, or even one which exists in $ZFC$! 

The example of 3.5 is clearly not first-countable.  In fact it can't be because of the following
result:

\proclaim{Theorem 3.6} If $X$ is a first-countable hereditarily Lindel\"of space with a small
diagonal, then $X$ has a $G_\d$-diagonal.  
\endproclaim
\demo{Proof} Suppose $X$ satsifies the hypotheses but not the conclusion.  Note that $X$
cannot have a countable $T_0$-separating open cover $\Cal U$, for otherwise
$$\Delta_X=\bigcap_{U\in \Cal U, n\in \o}[U^2 \cup(X\sm U_n)^2]$$ where $\{U_n:n\in \o\}$   
is a countable closed cover of $U$.  

Let $B(x,n)$, $n<\o$, be a countable base at $x$.  Then we can construct doubletons
$\{x_\a^0,x_\a^1\}$, $\a<\o_1$, such that $\{x_\a^0,x_\a^1\}$ is not separated by any member of
$\{B(x_\b^e,n):\b<\a, n<\o, e<2\}$. By Proposition 1.2(c) there is an open set $V$ with
$W=\{\a:x_\a^0\in V, x_\a^1\not\in V\}$ uncountable. For each $\a\in W$, choose $n_\a\in \o$
with $B(x_\a^0,n_\a)\subset V$.  There is $\d<\o_1$ such that $\{B(x_\b^0,n_\a):\b\in W\cap
\d\}$ covers $\{x_\a^0:\a\in W\}$.  But now if $\a\in W\sm\d$, then $\{x_\a^0,x_\a^1\}$ is
separated by $B(x_\b^0,n_\a)$ for some $\b\in W\cap \d$, contradiction. \qed
\enddemo

{\bf Remark.} We don't know if the above result remains true with ``hereditarily Lindel\"of"
weakened to ``Lindel\"of".  Zhou[Z] showed that the answer is ``yes"
under $CH$.  Also, Bennett and Lutzer [BL] showed in $ZFC$ 
that any Lindel\"of space with a small diagonal which is a subspace of
a linearly ordered space must have a $G_\d$-diagonal.
 
\head 4.  Countably compact spaces \endhead
Zhou mentions in [Z] that it is unknown if countably compact spaces
with a small diagonal must be metrizable (equivalently, have a
$G_\d$-diagonal), and he obtains some partial results related to
this question.  Shakhmatov [Sh] asks if one can at least show that 
they must be compact.   
In this section, we show that the statement ``Countably compact
spaces with a small diagonal are metrizable" is consistent with and
independent of $ZFC$. 

Let us note that the space $\o_1$ of countable ordinals does not
have a small diagonal, for the sequence $\<(\a,\a+1)\>_{\a<\o_1}$
in $\o_1^2$ converges to the diagonal.  More generally, Bennett and
Lutzer [BL] have shown that countably compact suborderable spaces
having a small diagonal are metrizable. 

The positive consistency result follows easily from the following
recent powerful result of Eisworth and Nyikos [EN]:

\proclaim{Theorem 4.1[EN]} The following statement $(*)$ is relatively 
consistent with \newline $ZFC+CH$: 

\noindent $(*)$\ \ A countably compact first-countable space is 
either compact or contains a copy of $\o_1.$  

\endproclaim

\proclaim{Theorem 4.2} If $CH$+$(*)$ holds, then countably compact spaces with a
small diagonal are metrizable.  \endproclaim
\demo{Proof}  Suppose $CH$+($*$) holds and $X$ is a countably compact space 
with a small diagonal which is not metrizable.  By CH and the Juh\'asz-
Szentmikl\'ossy result, $X$ is not compact.

{\it Case 1.  $X$ has a separable closed non-compact subspace $Y$.}
By CH, $Y$ has character not greater than $\o_1$.  Suppose some point
$y\in Y$ has character exactly $\o_1$.  The point $y$ cannot be a 
$G_\d$-point of $Y$, for if $U_n$, $n<\o$, is a sequence of neighborhoods
of $y$ with $\{y\}=\bigcap_{n<\o}U_n$ and $\cl{U_{n+1}}\subset U_n$ for all 
$n$, then it follows from countable compactness that $\{U_n\}_{n<\o}$ is a 
(countable) base at $y$.  So, if now $V_\a$, $\a<\o_1$, is a base at $y$, 
we can choose $y_\a\in \bigcap_{\b<\a}V_\a$, $y_\a\neq y$, and then 
$\<y_\a\>_{\a<\o_1}$ is a convergent $\o_1$-sequence, contradicting the 
small diagonal property.     

Thus $Y$ must be first-countable.  Since $Y$ is not compact, by ($*$) we
have that $Y$ contains a copy of $\o_1$.  But $\o_1$ has no small diagonal, 
contradiction.

{\it Case 2.  Every separable closed subset of $X$ is compact.}
Let $\Cal U$ be an open cover of $X$ with no finite (hence countable) subcover. 
Using the hypothesis of Case 2, one can easily construct points $y_\a$, $\a<\o_1$, and finite
subcollections $\Cal U_\a$ of $\Cal U$, such that $\Cal U_\a$ covers $\cl{\{y_\b:\b<\a\}}$ and
$y_\a\not\in \bigcup_{\b\leq \a}\cup \Cal U_\b$.  Let $Y=\bigcup_{\a<\o_1}\cl{\{y_\b:\b<\a\}}$. 
Then
$Y$ is countably compact, non-compact ($\bigcup_{\a<\o_1}\Cal U_\a$ is an open cover with no
finite subcover), and each $\cl{\{y_\b:\b<\a\}}$ is a compact, hence metrizable 
(by $CH$), open subspace.  It
follows
that $Y$ is first-countable.  Now ($*$) implies that $Y$ contains a copy of $\o_1$,
contradiction. \qed  
\enddemo

On the other hand, there are in some models countably compact spaces
with a small diagonal which are not metrizable.  Recall that a space
is {\it initially $\o_1$-compact} if every open cover of cardinality 
$\o_1$ or less has a finite subcover; equivalently, every subset of 
cardinality $\o$ or $\o_1$ has a complete accumulation point.

\proclaim{Theorem 4.3} Suppose $2^{\o_1}=\o_2$ and 
there is a subset $X$ of $2^{\o_2}$ of cardinality $\o_2$ satisfying:
\roster
\item"(a)" For every infinite subset $Y$ of $X$, there is $\g<\o_2$ such 
that $\{y\r (\o_2\sm \g): y\in Y\}$ is dense in $2^{\o_2\sm \g}$;
\item"(b)" For every $\o_1$-sequence $\<x_\a^0,x_\a^1\>_{\a<\o_1}$
in $X^2\sm\Delta$, and for every $\gamma<\o_2$, there is $\a<\o_1$ and
$n<\o$ with $x_\a^0(\g+n)\neq x_\a^1(\g+n)$.
\endroster
Then there is an initially $\o_1$-compact non-compact space with a small diagonal
but no $G_\d$-diagonal.  
\endproclaim

{\bf Remark.} A set $X$ satisfies the hypotheses of the above theorem, if,
e.g., it is both $HFD$ (to get (a)) and $HFD^2_\o$ (to get (b)) 
in $2^{\o_2}$.  Such a set, 
along with $2^{\o_1}=\o_2$, can be obtained by adding 
$\o_2$-many Cohen reals to a model of $2^{\o_1}=\o_2$.  
See [J$_2$] or [J$_3$] for information on $HFD$'s.  

\demo{Proof}  Let $\{x_\a:\a<\o_2\}$ index $X$, and let
$\{g_\a:\a<\o_2\}$ index all functions from ordinals less 
than $\o_2$ into $2$ (by $2^{\o_1}=\o_2$, there are not more than 
$\o_2$ of them).  Then define $z_\a\in 2^{\o_2}$ by 
$z_\a\r \operatorname{dom}(g_\a)=g_\a$ and $z_\a\r(\o_2\sm \operatorname{dom}(g_\a))=
x_\a\r(\o_2\sm \operatorname{dom}(g_\a))$.  We claim that the subspace 
$Z=\{z_\a:\a<\o_2\}$ of $2^{\o_2}$ is the 
desired example.  

It is easy to see that any infinite subset of $2^{\o_2}$ satisfying (a), 
which is satisfied by $Z$  since 
each $z_\a$ is the same as $x_\a$ beyond some ordinal $<\o_2$,
cannot be compact.
Initial $\o_1$-compactness of $Z$ also follows from (a).  The
proof of this is the same as the proof of Theorem 5.4 of [BG], which in
turn is a mild generalization of the proof due to Juh\'asz (see [J$_2$]) that  
an analogous construction in $2^{\o_1}$ yields a countably compact
space.  For the benefit of the reader, we repeat  the argument here. 

We need to show that every 
 subset of $Z$ of
cardinality $\o$ or $\o_1$ has a complete accumulation point.  
To this end,  let $\k\in \{\o,\o_1\}$, and suppose $Y\in [Z]^\k$.  
Using property (a),
it is easy to see that for each $W\in [Z]^\k$ there is $\d_W<\o_2$ such
that $\{w\r( \o_2\sm\d_W):z\in Z\}$ is $\k$-dense in $2^{\o_2\sm \d_W}$ (consider
splitting $W$ into $\k$-many infinite disjoint pieces).

Now we can find $\g<\o_2$ satisfying:
\roster
\item"(i)" The projection $\pi_\g:Y\to 2^\g$ is one-to-one;
\item"(ii)" For each $\s\in \operatorname{Fn}(\g,2)$, if $|Y\cap[\s]|=\k$ then
$\d_{Y\cap[\s]}<\g$.
\endroster
($\operatorname{Fn}(\g,2)$ is the set of all finite functions from a subset of $\g$
into $2$, and $[\s]=\{x\in 2^{\o_2}: x\text{ extends }\s\}$.)

Choose a complete accumulation point $g\in 2^\g$ of $\pi_\g(Y)$, and
let $z\in Z$ be an extension of $g$.  We claim that $z$ is a complete
accumulation point of $Y$.  Suppose $\s\in \operatorname{Fn}(\o_2,2)$ with $z\in [\s]$.
Let $W=Y\cap [\s\r \g]$; then $|W|=\k$.  By (ii), $\d_W<\g$, so $|W\cap
[\s\r(\o_2\sm\g)]|=\k$.  Thus $|Y\cap [\s]|=\k$. Thus $Z$ is initally
$\o_1$-compact.

Finally, that $Z$ has a small diagonal follows from (b), which is satisfied by $Z$
if ``for every $\g<\o_2$" is replaced by ``for sufficiently large
$\g<\o_2$".  For, suppose $\{z_\a^0,z_\a^1\}_{\a<\o_1}$ is an
 $\o_1$-sequence of doubletons of $Z$.  Write $\o_1$ as the union of 
$\o_1$-many disjoint uncountable sets $W_\a$, $\a<\o_1$.  For each $\a$,      
there is $\g_\a<\o_2$  such that the condition of (b) for the sequence 
$\<z_\b^0,z_\b^1\>_{\b\in W}$ is satisfied for
all $\g>\g_\a$.  Choose $\d<\o_2$, $\d> sup\{\g_\a:\a<\o_1\}$.  Then for 
each $\a<\o_1$, there is some $n_\a<\o$ and some $\b_\a\in W_\a$ with   
$z_{\b_\a}^0(\d+n_\a)\neq z_{\b_\a}^1(\d+n_\a)$.  Then for some $n<\o$,  
uncountably many doubletons of the original sequence differ on $\d+n$.  
Thus $Z$ has a small diagonal by Proposition 1.2.
\qed
\enddemo

We conclude this section with some questions:

I. Does $CH$ imply that countably compact spaces having a small
diagonal are metrizable?  What about $PFA$ or the statement ($*$) 
of Theorem 3.1?

II. Can there exist a non-metric countably compact space with a small diagonal
which is countably tight, or first-countable?  Dow[D] has
shown initially $\o_1$-compact, countably tight spaces with a small
diagonal are metrizable in models obtained by Cohen forcing over a model 
of $CH$.  The example of Theorem 4.3,
which does exist in some Cohen models, 
is not countably tight.  

Since Ostaszewski spaces (i.e., countably compact, locally compact, locally
\newline countable spaces in which every closed subset is either countable or 
co-countable) are particularly interesting examples of
first-countable countably compact non-compact spaces, we ask:

III.  Can there be an Ostaszewski space with a small diagonal?    

IV.  Can there be a first-countable perfect pre-image of $\o_1$  
with a small diagonal?  

The answer to Question III is ``no" under $MA_{\o_1}$, which 
kills Ostaszewskii spaces, and the answer to Question IV 
is ``no" under $PFA$, 
for then such a space would have to contain a copy of $\o_1$ [F$_2$].
Of course, by our Theorem 4.1 the answer to II,III, and IV is ``no"
under $CH$ together with the statement ($*$).

{\bf Remark.} O. Pavlov [P] recently obtained a positive solution 
to Question IV, namely, that there is such an example assuming 
axiom $\diamondsuit^+$.  This also gives a positive answer to Question II,
and answers in the negative the 
part of Question I about $CH$.  He has also shown 
that on the other hand 
there is no  finite-to-one perfect preimage of $\o_1$ with 
a small diagonal.

\head 5.  Locally compact spaces \endhead

Must locally compact spaces with a small diagonal have a
$G_\d$-diagonal?  Bennett and Lutzer [BL] showed that the 
answer is ``yes" for linearly ordered spaces (they actually use
a stronger assumption than small diagonal, 
but it turns out a little tweaking of the
argument gets it for small diagonal).  
Also, Zhou \cite{Z} obtained an example under $MA_{\o_1}$ showing 
that the answer can be ``no".   The purpose of this section is 
to show that the answer is ``no" in $ZFC$, i.e., we construct 
in $ZFC$ a locally
compact space with a small diagonal but no $G_\d$-diagonal.  

Both Zhou's example and ours are locally countable; their existence
depends essentially on the existence of almost disjoint families of
countable sets having certain combinatorial properties.  This makes
the problem for locally compact quite different than for compact;
but the combinatorics involved are natural and perhaps have some
interest in their own right.  

If $\Cal A$ is a collection of sets, let us say that a set $B$ is
{\it orthogonal} to $\Cal A$ if $B\cap A$ is finite for every $A\in
\Cal A$.  

\proclaim{Lemma 5.1} There is a locally countable, locally compact
$T_2$-space $X$ with a small diagonal but no $G_\d$-diagonal if
there is an almost disjoint family $\Cal A\subset [\k]^\o$ for some
cardinal $\k$ satisfying:
\roster
\item"(a)" For every uncountable subset $H$ of $\k$, there is an
uncountable subset  $H'$ of $H$ orthogonal to $\Cal A$;
\item"(b)" $\k$ is not the union of countably many subsets
orthogonal to $\Cal A$.
\endroster
\endproclaim

\demo{Proof}  Let $X=(\k\times 2)\cup \Cal A$, where $\k\times 2$
is a set of isolated points, and a neighborhood of $A\in \Cal A$ is
$\{A\}$ together with a cofinite subset of $A\times 2$.  

Let us see that $X$ has no $G_\d$-diagonal.  Suppose $\Cal G_n$,
$n<\o$, is a sequence of open covers of $X$.  For each $n$, let
$$B_n=\{\a\in \k: \forall G\in \Cal G_n(\{\<\a,0\>,\<\a,1\>\}
\not\subseteq G)\}.$$   Since for each $A\in \Cal A$, $\Cal G_n$
has an element $G$  containing all but finitely many points of
$A\times 2$, it follows $B_n$ is orthogonal to $\Cal A$.  Hence
there is $\a\in \k\sm \bigcup_{n\in \o}B_n$.  Then for each $n$,
$\<\a,1\>\in st(\<\a,0\>,\Cal G_n)$, whence $\Cal G_n$ cannot be a
$G_\d$-diagonal sequence for $X$.  

Now let us see that $X$ has a small diagonal.  Suppose
$\{\<x_\a,y_\a\>\}_{\a<\o_1}$ is an uncountable subset of $X^2\sm
\Delta$.  W.l.o.g., the $x_\a$'s are distinct.  

Suppose uncountably many $x_\a$'s are in $\Cal A$.  For each $A\in
\Cal A$, let $N(A)$ be a neighborhood of  $A$ not containing $y_\a$
if $x_\a=A$.     Then $$U=(\cup\{N(A)^2:A\in \Cal
A\})\cup\{\<z,z\>:z\in \k\times 2\}$$ is an open neighborhood of
$\Delta$ missing the uncountably many points with $x_\a\in \Cal A$. 

It remains to consider the case where uncountably many $x_\a$'s are
in $\k\times 2$.  Let $H=\{\g\in \k: \exists \a<\o_1\exists e<2
(x_\a=\<\g,e\>\}$.  $H$ contains an uncountable subset $H'$
orthogonal to $\Cal A$.  For each $A\in \Cal A$, let
$N'(A)=\{A\}\cup ((A\sm H')\times 2)$.  Let $U'$ be the
neighborhood of $\Delta$ defined as in the previous paragraph using
$N'(A)$ instead of $N(A)$.  Then $U'$ misses the uncountably many
points with $x_\a\in H'\times 2$.  That completes the proof.  \qed
    
\enddemo

\proclaim{Theorem 5.2}  There is an almost disjoint collection $\Cal A$
of countable subsets of $\o_1$ satisfying the conditions of Lemma
5.1, and hence there is a locally compact locally countable $T_2$-
space with a small diagonal but no $G_\d$-diagonal.  
\endproclaim

\demo{Proof}  For each limit ordinal $\a$ in $\o_1$, let $y_\a\in
\o_1^\o$ be such that $y_\a(n)$, $n\in \o$, is an increasing
sequence of ordinals with supremum $\a$.  Let $Y=\{y_\a:\a<\o_1\}$. 
Viewed as a subset of the metric space $\o_1^\o$, where $\o_1$ is
given the discrete topology, the space $Y$ was considered by A.H.
Stone \cite{St} in his non-separable Borel theory.  A pertinent fact
here is that every separable subspace of $Y$ is countable, so in
particular $Y$ contains no copy of a Cantor set.

Let $\Cal A\subset [Y]^\o$ be a maximal almost disjoint family of
sets $A$ that have a unique limit point $\o_1^\o\sm Y$.  We will
see that this $\Cal A$ satifies conditions (a) and (b) of Lemma 5.1.

Let $H\in [Y]^{\o_1}$.  In the metric topology, $H$ is not
separable, hence has an uncountable discrete (relative to $H$)
subset $H'$.  Since the set $K$ of limit points of $H'$ in 
$\o_1^\o$ is closed in that metric space, $K$ is $G_\delta$ and
hence $H'$ has an uncountable subset $H''$ which is closed in
$\o_1^\o$.  Clearly $H''$ is orthogonal to $\Cal A$.  Thus $\Cal A$
satisfies condition (a).

We check that $\Cal A$ satisfies condition (b).   Let $B\subset Y$
be orthogonal to $\Cal A$.  Condition (b) will follow if we show
that $S=\{\a:y_\a\in B\}$ is non-stationary.   Suppose $S$ is 
stationary; then by Stone's result \cite{St}, the closure of $B$ in
$\o_1^\o$  contains a copy of a Cantor set (in fact a copy of
$\o_1^\o$).  $Y$ does not contain a Cantor set, so 
some sequence $\{b_n:n\in\o\}$
of points of $B$ converges to some point of $\o_1^\o\sm Y$.  But this 
sequence must meet some member of $\Cal A$ in an infinite set, 
contradicting $B$ orthogonal to $\Cal A$.   \qed

\enddemo

{\bf Remark.}  S. Todor\v cevi\'c  independently discovered a (different)
almost disjoint collection $\Cal A$ satisfying the conditions of
Lemma 5.1.

\enddocument